**Divisor and Totient Functions Estimates**
**N. A. Carella, June, 2008.**

**Abstract:** New unconditional estimates of the divisor and totient functions are contributed to the literature. These results are consistent with the Riemann hypothesis and seem to solve the Nicolas inequality for all sufficiently large integers.

# 1 Introduction

The divisor function $\sigma(N) = \sum_{d \mid N} d$ and Euler function $\varphi(N) = \#\{ n < N : \gcd(n, N) = 1 \}$ are ubiquitous in number theory. Their product form representations

(i) $\sigma(N) = \prod_{p^{\alpha} \parallel N} (1 + p^2 + \Lambda + p^{\alpha})$ and (ii) $\varphi(N) = N \prod_{p \mid N} (1 - 1/p)$

unearth their intrinsic link to the distribution of prime numbers.

The divisor function is an oscillatory function, its value oscillates from its minimum $\sigma(N) = N + 1$ at prime integers $N$ to its maximum $\sigma(N) = c_0 N \log\log N$, some constant $c_0 > 1$, at extremely abundant integers $N$. An extremely abundant integer $N$ is an integer with lots of small prime factors and certain multiplicative structure. Similarly, the totient function $\varphi(N)$ is an oscillatory function, its value oscillates from its maximum $\sigma(N) = N - 1$ at prime integers $N$ to its minimum $\sigma(N) = N/c_0 \log\log N$, some constant $c_0 > 1$.

Currently the best unconditional estimates are the following:

***Theorem* 1.** ([R]) Let $N \in \mathbb{N}$, then $\sigma(N) < e^{\gamma} N \log\log N + \dfrac{.6482 N}{\log\log N}$ for $N \geq 3$.

***Theorem* 2.** ([RS]) Let $N \in \mathbb{N}$, then $N / \varphi(N) < e^{\gamma} N \log\log N + \dfrac{5}{2\log\log N}$ with one exception for $N = 2 \cdot 3 \cdots 23$.

On the other hand there are several conditional criteria; some of these are listed below. These results herein are stated in the notations of readily available sources such as [L], [S], et cetera, and other freely available papers.

**Theorem 3.** ([R])  (i) If the Riemann Hypothesis is true then for each $N \geq 5041$, $\sigma(N) < e^\gamma N \log \log N$ .

(ii) If the Riemann Hypothesis is false then there exists constants $0 < \beta < 1/2$ and $c > 0$ such that

$$\sum_{d \mid N} d \geq e^\gamma N \log \log N + \frac{cN \log \log N}{(\log N)^\beta} \tag{1}$$

holds for infinitely many $N$.

The parameter $1 - b < \beta < 1/2$ arises from the possibility of a zero $\rho \in \mathbb{C}$ of the zeta function $\zeta(s)$ on the half plane $b = \mathrm{Re}(\rho) > 1/2$. This in turns implies the existence of more primes per interval than expected if the Riemann hypothesis is valid. The effect of the zeros of the zeta function on the distribution of primes is readily revealed by the explicit formula.

**Theorem 4.**  ([L])  Let $H_N = \sum_{n=1}^{N} 1/n$ be the harmonic series. For each $N \geq 1$, the inequality $\sigma(N) < H_N + \exp(H_N) \log H_N$ is equivalent to the Riemann Hypothesis.

**Theorem 5.**  ([S])  The Riemann Hypothesis holds if and only if for all integers $N$ divisible by the fifth power of some prime, the inequality $\sigma(N) < e^\gamma N \log \log N$ holds.

**Theorem 6.**  ([NS])  Let $N_k = 2 \cdot 3 \cdot 5 \cdots p_k$ be the product of the first $k$ primes.

(i) If the Riemann Hypothesis is true then $N_k / \varphi(N_k) > e^\gamma \log \log N_k$ for all $k \geq 1$.

(ii) If the Riemann Hypothesis is false then both $N_k / \varphi(N_k) < e^\gamma \log \log N_k$ and $N_k / \varphi(N_k) > e^\gamma \log \log N_k$ occur for infinitely many $k \geq 1$.

Earlier works on this topic include the works of Ramanujan and other on abundant numbers, see [RJ], [AE], and recent related works appeared in [BR], [BS], and [WZ]. The new contributions to the literature are the unconditional estimates stated below.

**Theorem 7.**  Let $N_k = 2 \cdot 3 \cdot 5 \cdots p_k$. Then $N_k / \varphi(N_k) > e^\gamma \log \log N_k$ for all sufficiently large integer $N_k$.

**Theorem 8.**  Let $N \in \mathbb{N}$, then $\sigma(N) < e^\gamma N \log \log N$ for all sufficiently large integers $N \geq 1$ such that $P(N) < \left(1 - 1/9 \log \log N\right) \log N$ .

**Theorem 9.**  If $N = 2^\nu M$ , where $M$ is odd and $2^\nu \leq (\log \log N)^2$, then $\sigma(N) < e^\gamma N \log \log N$ for all such $N$ but finitely many exceptions.

These unconditional results are consistent with the Riemann hypothesis, and seem to prove the Nicolas inequality, Theorem 6-i, for all sufficiently large integers. Just a finite number of cases remain unresolved as possible counterexamples. The proofs of Theorems 7, 8 and 9 are given in Sections 10, 11, and 12. The other sections are background and supplemental materials focusing on various characteristic of the divisor function, totient function, and other arithmetic functions.

## 2 Properties of the Divisor Function and Identities

A more general version of the divisor function takes the form $\sigma_s(N) = \sum_{d\,|\,N} d^s$, where $s \in \mathbb{C}$ is a complex number. The most common cases are the divisor function $\sigma_0(N) = \sum_{d\,|\,N} 1$, and the sum of divisor function $\sigma(N) = \sigma_1(N) = \sum_{d\,|\,N} d$ .

### Properties of the Divisor Function

Let $N = p_1^{\alpha_1} \cdot p_2^{\alpha_2} \Lambda \ p_k^{\alpha_k}$, and let the symbol $p^\alpha \,\|\, N$ denotes the maximum prime power divisor.

1) $\sigma_0(N) = (\alpha_1 + 1)(\alpha_2 + 1)\Lambda \ (\alpha_k + 1)$ .

2) $\sigma_s(N) = \prod_{p^\alpha\,\|\,N} \sigma_s(p^\alpha) = \prod_{p^\alpha\,\|\,N} \dfrac{p^{(\alpha+1)s} - 1}{p - 1}$, \qquad Multiplicative.

3) $\sigma_s(N) = \sum_{d\,|\,N} d^{-s}, \quad N^s = \sum_{d\,|\,N} \mu(N/d)\sigma_s(d)$, \qquad Mobius inversion pair.

4) $\sigma(N)^2 = N \sum_{d\,|\,N} d^{-1}\sigma(d^2)$ .

5) If $d \geq 1$ divides $N$ then $\sigma(d)/d \leq \sigma(N)/N$ .

6) $\sum_{d\,|\,N} \mu(d)\sigma(d) = (-1)^k \, p_1 p_2 \Lambda \ p_k$ .

7) If $\gcd(M, N) > 1$, then $\sigma_s(MN) < \sigma_s(M)\sigma_s(N)$ .

8) If $p\,|\,N$ and $s > 0$, then $\sigma_s(N) = \sigma_s(N)\sigma_s(p) - p^s\sigma_s(N/p)$, \qquad Two terms recursive formula.

9) $\sigma_s(N)^2 = \sum_{d\,|\,N} d^s \sigma_s(N^2/d^2)$ and $\sigma_s(M)\sigma_s(N) = \sum_{d\,|\,\gcd(M,N)} d^s \sigma_s(MN/d^2)$ .

10) $\sigma_s(N) < \prod_{p\,|\,N} p^s/(p^s - 1)$, and $\sigma_s(N) < \zeta(s)N^s$ .

***Theorem* 10.** For every $s \in \mathbb{C}$ and integer $N \geq 1$, the minimum order of $\sigma_s(N)$ is $N^s$. For $\mathrm{Re}(s) > 1$, the maximum order is $\zeta(s)N^s$. For $s = 1$, the maximum order is $e^\gamma N \log\log N$ . For $s < 1$, the maximum order is

$$\sigma_s(N) \leq N^s e^{(1+o(1))(\log)^{1-s}/(1-s)\log\log N}, \qquad (2)$$

and the opposite inequality holds for infinitely many integers.

The case $\mathrm{Re}(s) > 1$ is easy to resolve by mean of the multiplicative formula $\sigma_s(N) = \prod_{p^\alpha\,\|\,N}(1 + p^{2s} + \Lambda + p^{\alpha s})$, and the case $s = 1$ is Gronwall's Theorem, more details on this appears in [TN, p. 88].

### Properties of the Totient Function

Let $N = p_1^{\alpha_1} \cdot p_2^{\alpha_2} \Lambda \ p_k^{\alpha_k}$, and let $\varphi_s(N) = \#\{\, (n_1,...,n_s) : \gcd(n_i, N) = 1, n_i < N \,\}$ for $s > 0$.

1) $\varphi(N) = N \prod_{p\,|\,N}\left(1 - 1/p\right)$ .

2) $\varphi(N) = N \sum_{d\,|\,N} \mu(d)d^{-1}$ and $\varphi(N) = \sum_{d\,|\,N} \varphi(d)$ .

3) $\varphi_s(N) = N^s \prod_{p\,|\,N}\left(1 - 1/p^s\right), \quad N^s = \sum_{d\,|\,N} \mu(N/d)\varphi_s(d)$ .

4) $\dfrac{N}{\varphi(N)} = \left(\dfrac{\mu(m)}{m}\right)^{-1} \displaystyle\sum_{d\,\mid\,N,\,m\mid d} \dfrac{\mu(d)^2}{\varphi(d)}$ for any $m \mid N$.

## 3 Representations of the Divisor and Totient Functions

Many of the important properties of the divisor and totient functions can be considered representations of these functions. All these representations are useful in the analysis of these functions. A few of them are recorded here.

**Proposition 11.** Let $N$ be an integer, and let the symbol $p^{\alpha} \parallel N$ denotes the maximum prime power divisor. Then

(i) $\sigma(N) = \dfrac{N}{\varphi(N)} \displaystyle\prod_{p^{\alpha} \parallel N}(1 - 1/p^{\alpha+1})$ .

(ii) $\sigma(N) = \sigma_2(N) \displaystyle\prod_{p^{\alpha} \parallel N} \dfrac{p+1}{p^{\alpha+1}+1}$ .

(iii) $\sigma(N) = \sigma_3(N) \displaystyle\prod_{p^{\alpha} \parallel N} \dfrac{p^2+p+1}{p^{2(\alpha+1)}+p^{\alpha+1}+1}$ .

The first and second of these representations are well known, and the third appear to be new. More general versions of these representations are possible, but are let to the reader to work out. Other related identities are given below.

**Proposition 12.** Let $N$ be an integer, and let the symbol $p^{\alpha} \parallel N$ denotes the maximum prime power divisor. Then

(i) $\dfrac{\sigma(N)}{N} = \dfrac{\sigma_2(N)}{N^2} \displaystyle\prod_{p^{\alpha} \parallel N} \left(\dfrac{p^{3\alpha+2}+p^{\alpha}}{p^2+1}\right)\left(\dfrac{p^3+p^2+p+1}{p^{3(\alpha+1)}+p^{2(\alpha+1)}+p^{\alpha+1}+1}\right)$,

(ii) $\dfrac{\sigma(N)}{N} = \dfrac{\sigma_3(N)}{N^3} \displaystyle\prod_{p^{\alpha} \parallel N} \left(\dfrac{p^{5\alpha+3}+p^{2\alpha}}{p^3+1}\right)\left(\dfrac{p^5+p^4+p^3+p^2+p+1}{p^{5(\alpha+1)}+p^{4(\alpha+1)}+p^{3(\alpha+1)}+p^{2(\alpha+1)}+p^{\alpha+1}+1}\right)$,

(iii) $\dfrac{\sigma(N)}{N} < \displaystyle\prod_{p\,\mid\,N} \left(1-\dfrac{1}{p^2}\right)^{-1}\left(1+\dfrac{1}{p}\right) < \zeta(2)\displaystyle\prod_{p\,\mid\,N}\left(1+\dfrac{1}{p}\right)$ .

(iv) $\dfrac{\sigma(N)}{N} < \displaystyle\prod_{p\,\mid\,N} \left(1-\dfrac{1}{p^3}\right)^{-1}\left(1+\dfrac{1}{p}+\dfrac{1}{p^2}\right) < \zeta(3)\displaystyle\prod_{p\,\mid\,N}\left(1+\dfrac{1}{p}+\dfrac{1}{p^2}\right)$ .

Proof: To verify (ii) consider computing the power divisor function $\sigma_6(N)$ in two ways:

$$\sigma_6(N) = \prod_{p^{\alpha} \parallel N}\left(\dfrac{p^{6(\alpha+1)}-1}{p^6-1}\right) = \sigma(N)\prod_{p^{\alpha} \parallel N}\left(\dfrac{p^{5(\alpha+1)}+p^{4(\alpha+1)}+p^{3(\alpha+1)}+p^{2(\alpha+1)}+p^{\alpha+1}+1}{p^5+p^4+p^3+p^2+p+1}\right),$$

and

$$\sigma_6(N) = \prod_{p^{\alpha} \parallel N}\left(\dfrac{p^{6(\alpha+1)}-1}{p^6-1}\right) = \sigma_3(N)\prod_{p^{\alpha} \parallel N}\left(\dfrac{p^{3(\alpha+1)}+1}{p^3+1}\right).$$

Now rearrange the products. And to verify (iv), use routine algebraic manipulations to simplify the product (ii), and then apply Property 10 in Section 2. ∎

There are several ways of establishing these results.

## 4 Related Arithmetic Functions and Average Orders

***Theorem 13.*** For $x \geq 1$, the average orders of the divisor functions are

(i) $\displaystyle\sum_{n \leq x} \sigma_0(n) = x \log x + (2\gamma - 1)x + O(x^\theta)$, where $\alpha < 1/2$.

(ii) $\displaystyle\sum_{n \leq x} \sigma_s(n) = \frac{\zeta(s+1)}{s+1} x^{s+1} + O(x^\beta)$, where $\beta = \max\{1, s\}$.

Confer [AP, p 61] for the analysis. The error terms of arithmetic functions are well studied problems in number theory, extensive details are given in [IV]. For example, the latest improvement on the error term $\Delta(n) = O(x^\theta)$ in the expression for the number of divisor appears to be $\theta = 23/73 + \varepsilon$, Huxley 1993.

The analysis and determination of the average order of an arithmetic function over subsequences of integers are significantly longer. Two recently determined cases over the sequence of binomial coefficients are included here, the reader should confer the paper for the proofs.

***Theorem 14.*** ([LA]) Let $N$ be an integer. Then

(i) $\displaystyle\frac{1}{N+1} \sum_{k=0}^{N} \sigma\binom{N}{k} \Big/ \binom{N}{k} = O(\log \log \log N)$,

(ii) $\displaystyle\frac{1}{N+1} \sum_{k=0}^{N} \varphi\binom{N}{k} \Big/ \binom{N}{k} = O(1/\log \log \log N)$.

***Theorem 15.*** For $x \geq 1$, the average orders of the totient function is $\displaystyle\sum_{n \leq x} \varphi(n) = (6/\pi^2)x^2 + O(x \log x)$.

## 5 A few Prime Numbers Results

The $n$th prime in the sequence of prime 2, 3, 5, 7, ... is denoted by $p_n$, and let $\omega(N) = \#\{p \mid N : p \text{ is prime}\}$ be the number of prime divisors counting function.

***Theorem 16.*** (Cipolla 1902) For $n \geq 2$, $n \log n \leq p_n \leq n(\log n + \log \log n - 1)$. In particular, the power series expansion is

$$p_n = n\left(\log n + \log \log n - 1 + \frac{\log \log n - 2}{\log n} - \Lambda\right).$$

Confer [DT] for recent developments in this area. The work of Hardy and Ramanujan on the function $\omega(N) = \#\{p \mid N\}$ culminated in the probabilistic result given below.

***Theorem 17.*** Let $\varepsilon > 0$, and $N \geq 1$. Then

(i) $(1-\varepsilon)\dfrac{\log N}{\log \log N} < \omega(N) < (1+\varepsilon)\dfrac{\log N}{\log \log N}$,

(ii) $\displaystyle\sum_{n \leq x} \omega(n) = x \log x + cx + O(x/\log x)$.

***Theorem 18.*** (Erdos-Kac 1940) The random variable $\omega(N)$ is normal with mean $\log \log N$ and standard deviation $\sqrt{\log \log N}$. Let $A(x) = \#\{n \leq x : \log \log n + a\sqrt{\log \log n} < \omega(n) < \log \log n + b\sqrt{\log \log n}\}$, where $a, b$ $\in \mathbb{R}$. Then $\dfrac{A(x)}{x} = \dfrac{1}{\sqrt{2\pi}} \displaystyle\int_a^b e^{-z^2/2} dz$ as $x$ tends to infinity.

The current perspective on the analysis of the function $\omega(N)$ is discussed in fine details in [GR].

The Chebychev functions are defined by

$$\vartheta(x) = \sum_{p \leq x} \log p \quad \text{and} \quad \psi(x) = \sum_{p \leq x} \alpha \log p,$$

where $p^\alpha < x$, but $p^{\alpha+1} > x$. The first is the logarithm of the product of all the primes $\leq x$, and the second is logarithm of the lowest common multiple of the integers $\leq x$. A related function appears in the explicit formula

$$\psi_0(x) = x - \sum_\rho \frac{x^\rho}{\rho} - \frac{\zeta'(0)}{\zeta(0)} - \frac{1}{2} \log\left(1 - \frac{1}{x^2}\right),$$

where $\rho = 1/2 + it$ are the zeroes of the zeta function $\zeta(s) = \sum_{n \geq 1} n^{-s}$, and $\psi_0(x) \approx \psi(x)$.

**Theorem 19.** (Chebyshev 1850) Let $x > 1$. The (i) $c_0 x < \vartheta(x) < c_1 x$ for some constants $c_0, c_1 > 0$.
(ii) $\psi(x) = \vartheta(x) + O(x^{1/2} \log x)$,          (iii) $\lim_{x \to \infty} x / \vartheta(x) = 1$.

The constants are approximately one, and the error term $E(x) = x - \vartheta(x)$ tends to infinity as $x$ tends to infinity. The function $\vartheta(x)$ is monotonically increasing, but $\psi(x)$ is an oscillating function of $x$. The peaks and valleys of the oscillation are known to satisfy $|\psi(x) - x| > c\sqrt{x}$ infinitely often, $c > 0$ constant.

**Lemma 20.** Let $x \geq 2$, then $\left| x - \vartheta(x) \right| \leq O(x \log^{-A} x)$ and $\left| x - \psi(x) \right| \leq O(x \log^{-A} x)$ some $A \geq 1$.

Proof: Confer [RS], [SC] and [DT] for other sharper estimates too.       ∎

**Proposition 21.** Let $N = 2^{v_1} \cdot 3^{v_2} \Lambda \ p_k^{v_k} \in \mathbb{N}$ and let $f(x) > 0$ be an strictly increasing function of $x$ such that $f(x) / x = o(x)$. Then the followings hold.
(i) If $\log N < p_k < (1 + 1/f(\log N)) \log N$ then $\log \log p_k - \log \log \log N = O(1/(f(\log N) \log \log N))$.
(ii) If $p_k < (1 - 1/f(\log N)) \log N$ then $\log \log p_k - \log \log \log N = O(1/(f(\log N) \log \log N))$.

Proof of (i): Write $p_k \leq (1 + 1/f(\log N)) \log N$, and simply compute the logarithmic difference:

$$\log\left(\frac{\log p_k}{\log \log N}\right) \geq c_0 \log\left(\frac{\log((1 + 1/f(N)) \log N)}{\log \log N}\right) = c_0 \log\left(\frac{\log \log N + \log(1 + 1/f(N))}{\log \log N}\right)$$

$$= c_0 \log\left(1 + \frac{\log(1 + 1/f(N))}{\log \log N}\right) \geq c_1 \frac{\log(1 + 1/f(N))}{\log \log N} \geq \frac{c_2}{f(\log N) \log \log N},$$

where $c_i >$ are constants. The proof of (ii) is similar.       ∎

The unconditional case uses the function $f(x) = c \log^A x$, some constants $A > 0$, $c > 0$, and the corresponding logarithm difference is

$$\log \log \log N - \log \log p_k = O(1/(\log \log N)^{A+1}).$$

In contrast, the conditional on the Riemann hypothesis case, uses the function $f(x) = cx^{1/2}$, and the corresponding logarithm difference is

$$\log \log \log N - \log \log p_k = O\left(1/(\log^{1/2} N \log \log N)\right).$$

## 6 Finite and Asymptotic Results

**Theorem 22.** For every integer $N \geq 3$, there is an absolute constant $c_0 > 0$ such that

(i) $\sigma(N) < c_0 N \log \log N$ 　　　　　　　　(ii) $\varphi(N) > N / c_0 \log \log N$.

Proof of (ii): Taking the logarithm and rearranging return

$$\log \frac{\varphi(N)}{N} = \sum_{p \mid N} \log\left(1 - \frac{1}{p}\right) = -\sum_{p \mid N} \frac{1}{p} + \sum_{p \mid N}\left(\log\left(1 - \frac{1}{p}\right) + \frac{1}{p}\right) \geq -\sum_{p \mid N} \frac{1}{p} - c_1.$$

The last step follows from

$$\sum_{p \mid N}\left(\log\left(1 - \frac{1}{p}\right) + \frac{1}{p}\right) = -\sum_{p \mid N}\sum_{n=2}^{\infty} \frac{1}{np^n} > -\frac{1}{2}\sum_{p \geq 2}\sum_{n=2}^{\infty} \frac{1}{p^n} = -\frac{1}{2}\sum_{p \geq 2} \frac{1}{p(p-1)} = -c_1.$$

Since $N \geq \prod_{p \mid N} p \geq \prod_{\log N \leq p \mid N} p \geq (\log N)^m$, where $m = \#\{p \mid N : p \geq \log N\}$, and

$$\sum_{p < \log N, p \mid N} \frac{1}{p} \leq \sum_{p < \log N} \frac{1}{p} = \log \log \log N + c_2 \quad \text{and} \quad \sum_{\log N \leq p \mid N} \frac{1}{p} \leq \frac{m}{\log N} \leq \frac{1}{\log \log N} \leq \frac{1}{\log \log 3}.$$

Combining these data give

$$\log \frac{\varphi(N)}{N} > -\log \log \log N - c_3.$$

Now taking the inverse logarithm yields the claim. 　　　　　　　　　　　■

The techniques just used above are quite standard in the literature, see [RD, p. 614], [SH, p. 341], [NT, p. 278], et cetera.

**Theorem 23.** (Landau 1903)　Let $N >$ be an integer. Then

(i) $\limsup\limits_{N \to \infty} \dfrac{N}{\varphi(N) \log \log N} = e^{\gamma}$, 　(ii) $\limsup\limits_{\text{odd } N \to \infty} \dfrac{N}{\varphi(N) \log \log N} = \dfrac{e^{\gamma}}{2}$, 　(iii) $\limsup\limits_{\text{squarefree } N \to \infty} \dfrac{N}{\varphi(N) \log \log N} = \dfrac{6e^{\gamma}}{\pi^2}$.

**Theorem 24.** (Gronwall 1913)　Let $N >$ be an integer. Then

(i) $\limsup\limits_{N \to \infty} \dfrac{\sigma(N)}{N \log \log N} = e^{\gamma}$, 　(ii) $\limsup\limits_{\text{odd } N \to \infty} \dfrac{\sigma(N)}{N \log \log N} = \dfrac{e^{\gamma}}{2}$, 　(iii) $\limsup\limits_{\text{squarefree } N \to \infty} \dfrac{\sigma(N)}{N \log \log N} = \dfrac{6e^{\gamma}}{\pi^2}$.

Another pair of related formulas are also recorded below.

**Theorem 25.** (Martens, 1874)  The following asymptotic formulas hold:

(i) $\displaystyle \lim_{n \to \infty} \frac{1}{\log n} \prod_{p \le p_n} \left(1 - 1/p\right)^{-1} = e^{\gamma}$,  (ii) $\displaystyle \lim_{n \to \infty} \frac{1}{\log n} \prod_{p \le p_n} \left(1 + 1/p\right) = \frac{6}{\pi^2} e^{\gamma}$.

A function $f(n)$ has *normal order* $F(n)$ if $f(n)$ is approximately $F(n)$ for almost all values of $n > 0$, confer [HW, p. 356].

**Theorem 26.**  The normal orders of $\sigma(N)$ and $\varphi(N)$ are

$$\sigma(N) = e^{\gamma} N \log\log\log N \quad \text{and} \quad \varphi(N) = N / e^{\gamma} \log\log\log N \text{ respectively.} \tag{3}$$

Proof: By the Erdos-Kac Theorem, the random variable $\dfrac{\omega(N) - \log\log N}{\sqrt{\log\log N}}$ is normal with mean $\mu = 0$ and standard deviation $\sigma = 1$. Therefore, $\omega(N) \le 1.5 \log\log N$ for almost all sufficiently large integers $N$, this follows from the symmetric version of the Chebychev's inequality. Next from the estimate $\omega(N) = \pi(p_n) = n \le 1.5 \log\log N$, it follows that the $n$th prime $p_n$ satisfies the inequality

$$p_n \le 2n \log n \le (3\log\log N)\log(1.5 \log\log N) \le 6 \log\log N. \tag{4}$$

Put $x = 6 \log\log N$ and apply Merten's formula (Theorem 31-i) to this data. Then

$$\frac{\varphi(N)}{N} = \prod_{p \mid N}(1 - 1/p) = \frac{e^{-\gamma}}{\log x}\left(1 + O(1/\log x)\right) = \frac{e^{-\gamma}}{\log\log N}\left(1 + O(1/\log\log\log N)\right) \tag{5}$$

for almost all sufficiently large integer $N$. The proof for $\sigma(N)$ is derived from this via the relation $\sigma(N)/N = N/\varphi(N)\prod_{p^{\alpha} \| N}(1 - 1/p^{\alpha+1})$.  ∎

The fact that $\sigma(N)$ is a normal random variable is easily seemed by means of the Duncan's formula below. More precisely, for squarefree integers $\sigma(N)$ is bounded by a linear transformation $\sigma(N)/N < a\omega(N) + b$ of the normal random variable $\omega(N)$, where $a, b > 0$ are constants.

**Proposition 27.**  ([DN]) Let $N \in \mathbb{N}$ be a squarefree integer. Then $\sigma(N) < \dfrac{\pi^2}{6} N\left(1 + \omega(N)\log 2\right)$.

Proof: Rewrite the sum of divisors function as $\sigma(N) = N \sum_{d \mid N} \dfrac{1}{d} \le N \sum_{d=1}^{\sigma_0(N)} \dfrac{1}{d} \le NH(\sigma_0(N))$. Proceed to replace the estimates the number of divisors $\sigma_0(N) \le 2^{2\log N / \log\log N}$ for any integer $N \in \mathbb{N}$ or $\sigma_0(N) = 2^{\omega(N)}$ for squarefree

integers, and the harmonic series $H(N) = \sum_{n=1}^{N} 1/n \le 1 + \log N$ to sharpen the estimate:

$\sigma(N) \le N\left(1 + \log\left(2^{2\log N / \log\log N}\right)\right) = N\left(1 + \log(4)\log N / \log\log N\right)$ for any integers. In particular,

$\sigma(N) \le N\left(1 + \log(2)\omega(N)\right)$ for squarefree integers. This proves the claim.  ∎

## 7 Highly Composite and Abundant Integers

An integer $N$ is called *deficient* if $\sigma(N) < 2N$, *perfect* if $\sigma(N) = 2N$, and *abundant* if $\sigma(N) > 2N$. More generally, the integer is called $m$-abundant or multiperfect if $\sigma(N) = mN$.

**Proposition 28.** If $N$ is abundant then $mN$ is $(m+1)$-abundant for all $m \geq 1$.

To see this, observe that $m$ has at least two divisors, so $\sigma(mN) \geq \sigma(N) + m\sigma(N) = (m+1)\sigma(N)$.

An integer is called highly composite if the relation $\sigma_0(M) < \sigma_0(N)$ holds for every integer $M < N$. Highly composite Integers $N = 2^{v_2} \cdot 3^{v_3} \cdot 5^{v_5} \Lambda \ p_n^{v_p}$ must have monotonically decreasing exponents $v_2 \geq v_3 \geq v_5 \geq \cdots \geq v_p \geq 1$, and every prime up to the $n$th prime $p_n$ included in the product.

Given any exponents vector $u_1, u_2, \ldots, u_n \geq 1$, there are infinitely many integers $N = q_1^{u_1} \cdot q_2^{u_2} \cdot q_3^{u_3} \Lambda \ q_n^{u_n}$ with the same number of divisors $\sigma_0(N) = (u_1+1)(u_2+1)\Lambda \ (u_n+1)$ since permuting the exponents vector and/or varying the primes $q_i$ does not change $\sigma_0(N)$. However, there is only one or which is highly composite as specified in the definition.

Let $v_1 \geq v_2 \geq \Lambda \geq v_k \geq 1$ be a fixed exponents vector. A *colossally abundant* integer $N = 2^{v_1} \cdot 3^{v_2} \cdot 5^{v_3} \Lambda \ p_k^{v_k}$ maximizes the divisor function on the infinite set of integers $S(v_1, \ldots v_k) = \left\{ M = q_1^{v_1} \cdot q_2^{v_2} \Lambda \ q_k^{v_k} : q_i \text{ prime} \right\}$. Specifically, $\sigma(N) > \sigma(M)$ for any $M \in S(v_1, \ldots, v_k)$ such that $M > N$, a detailed discussion appears in [S]. A related definition states that a colossally abundant number is an integer for which

$$\frac{\sigma(N)}{N^{1+\varepsilon}} \geq \frac{\sigma(M)}{M^{1+\varepsilon}} \text{ for all } 1 \leq M < N, \text{ and } \varepsilon > 0.$$

The multiplicative structure of a colossally abundant number is specified by $N = 2^{v_2(\varepsilon)} \cdot 3^{v_3(\varepsilon)} \cdot 5^{v_5(\varepsilon)} \Lambda \ p_n^{v_p(\varepsilon)}$, where the exponents satisfy $v_2(\varepsilon) \geq v_3(\varepsilon) \geq \Lambda \geq v_p(\varepsilon) \geq 1$ and $v_p(\varepsilon) = \left[\log_p(p^{1+\varepsilon}-1)/(\log_p(p^{\varepsilon}-1))\right] - 1$. Numerical results and an algorithm for generating colossally abundant integers is described in [BS].

# 8 Harmonic and Quasiharmonic series
The *harmonic series* is the sum of the reciprocals of the positive integers up to $x$. A *quasiharmonic series* is a sum of the reciprocals of $x$ numbers less than or equal to $x$. The summation techniques used to estimate these series are mostly elementary. On the other hand the applications of these series in the mathematical sciences run deep.

A selection of useful series are recorded here, the proofs are scattered in the literature but easily available.

**Theorem 29.** For every positive number $x \geq 286$, the followings hold:

(i) $\sum_{n \leq x} \frac{1}{n} = \log x + \gamma + O(1/x)$.

(ii) $\sum_{n \equiv a \bmod q \text{ and } n \leq x} \frac{1}{n} = \frac{\log x}{\varphi(q)} + \gamma_{a,q} + O(1/x)$.

(iii) $\sum_{n \leq x \text{ and } \gcd(n,N)=1} \frac{1}{n} = \frac{\varphi(N)}{N} \log x + o(\frac{\varphi(N)}{N} \log x)$.

(iii) $\sum_{\text{squarefree } n \leq x} \frac{1}{n} = \frac{6}{\pi^2} \log x + c + O(x^{-1/2})$,

where $\gamma$ and $\gamma_{a,q}$ are the Euler's constants, see Section 15.

**Theorem 30.** For $x \geq 2$ and $k \geq 0$,

(i) $\sum_{n \le x} \frac{\log^k n}{n} = \frac{\log^{k+1} x}{k+1} + c + O((\log^k x)/x)$.     (ii) $\sum_{2 \le n \le x} \frac{1}{n \log n} = \log \log x + c + O(1/x \log x)$.

These are proved using summation methods, such as the Euler-Maclaurin Formula, see [RM, p. 234] and similar literature.

**Theorem 31.**   (Mertens 1874) Let $x \ge 1$. Then

(i) $\sum_{p \le x} \frac{1}{p} = \log \log x + B + R(x)$,     (ii) $\sum_{p \le x, \, p \equiv a \bmod q} \frac{1}{p} = \frac{1}{\varphi(q)} \log \log x + B_{a,q} + O_q(1/\log x)$,

$B = .2614972128\ldots$, and $B_{a,q}$ is Mertens constant on arithmetic progression, Section 15.

These are recent versions of the original work. The error term has the form $R(x) = O(1/\log^B x)$ or better for some constant $B > 0$. Currently the best estimate are

$$R(x) = \pm \left( \frac{1}{10 \log^2 x} + \frac{4}{15 \log^2 x} \right) \text{ or } R(x) = \pm \frac{3 \log x + 4}{8\pi \sqrt{x}} \tag{6}$$

unconditionally or modulo the Riemann hypothesis respectively, see [VL]. There are various analytical formula for the constant $B = \gamma - \sum_{p \ge 2} \sum_{n=2}^{\infty} \frac{1}{np^n}$, see [HW, p. 351], [VL, p. 3]. About four different proofs of Theorem 31-i are given in [DN, p. 56], [NK, p 126] and [RS, p. 68]. A new proof is given here, it do away with the asymptotic term and it is far simpler and very practical in applications.

**Theorem 32.**   Let $x \ge 3$. Then $\sum_{p \le x} 1/p \le \log \log x - \log \log 2 + 1/2 \log 2$.

Proof:  Let $x = p_k$ be the $k$th prime number. By Theorem 16, the $n$th prime satisfies the inequality $n \log n \le p_n \le n \log(n \log n)$ for $n \ge 2$. Summing over all the reciprocal primes returns

$$\sum_{n=2}^{k} \frac{1}{n \log(n \log n)} \le \sum_{p \le p_k} \frac{1}{p} \le \sum_{n=2}^{k} \frac{1}{n \log n}$$

for all $k \ge 4$. Now by the integral test, it follows that

$$\sum_{n=2}^{x} \frac{1}{n \log(n \log n)} \le \sum_{p \le x} \frac{1}{p} \le \sum_{n=2}^{x} \frac{1}{n \log n} \le \frac{1}{2 \log 2} + \int_{2}^{x} \frac{1}{x \log x} dx \le \log \log x - \log \log 2 + 1/2 \log 2. \quad \blacksquare$$

The left side of this estimate can be evaluated in terms of logarithm and the exponential integral, but it is let as an exercise for the reader to improve the analysis.

**Theorem 33.**   Let $x \ge 2$. Then

(i) $\sum_{p \le x} \frac{1}{p-1} = \log \log x + A_- + R(x)$,     (ii) $\sum_{p \le x} \frac{1}{p+1} = \log \log x + A_+ + R(x)$,

where $A_\mu = B + \sum_p 1/p(p \mu 1)$ is a constant.

Proof: Rewrite the finite sums as $\sum_{p \leq x} 1/(p+a) = \sum_{p \leq x} 1/p + \sum_{p \leq x} 1/p(p+a)$, and then apply Mertens formula, Theorem 31-i. ∎

**Proposition 34.** Let $x \geq 2$, then $\sum_{p > x} \sum_{n=2}^{\infty} \dfrac{1}{np^n} = O(1/x \log x)$.

Proof: Put $x = k \log k$, and use the fact that the $n$th prime $p_n$ satisfies $an \log n \leq p_n \leq bn \log n$, some constants $a$, $b > 0$. Then

$$\sum_{p > x} \sum_{n=2}^{\infty} \frac{(-1)^n}{np^n} = \sum_{n=2}^{\infty} \sum_{p > x} \frac{1}{np^n} = \sum_{n=k+1}^{\infty} \frac{1}{2p_n^2} + \sum_{n=k+1}^{\infty} \frac{1}{3p_{n+1}^3} + \Lambda = O(1/(k \log^2 k)).$$

This is obtained using the integral approximation

$$\sum_{n=k+1}^{\infty} \frac{1}{p_n^2} \leq \sum_{n=k+1}^{\infty} \frac{1}{(n \log n)^2} \leq \frac{1}{(k \log k)^2} + \int_k^{\infty} \frac{dt}{(t \ln t)^2} = O(1/(k \log^2 k)).$$
∎

## 9 Harmonic and Quasiharmonic Products

The harmonic product $\prod_{p \leq x}(1 - 1/p)$ is ubiquitous in number theory. It has a natural link to the Euler function, sieve of Eratosthenes, and other related concepts. In the case of the sieve of Eratosthenes, all multiples of the primes $p \leq x^{1/2}$ are removed from the set of integers $\{1, 2, 3, \ldots, x\}$. The remaining subset of primes $\{p \leq x\}$ contains approximately $x \cdot \prod_{p \leq x}(1 - 1/p)$ primes. Similarly, in the case of the Euler function calculation, the nonharmonic product $\prod_{p \mid N}(1 - 1/p)$ is the probability that an integer $m < N$ is relatively prime to $N$. Accordingly $N \cdot \prod_{p \mid N}(1 - 1/p)$ is exactly the number of integers $m < N$ relatively prime to $N$.

Accurate estimates of this product and related products are essential in a variety of calculations in number theory.

**Theorem 35.** (Mertens 1874) Let $x \geq 286$. Then

(i) $\prod_{p \leq x}(1 - 1/p) > \dfrac{e^{-\gamma}}{\log x}\left(1 - \dfrac{1}{2(\log x)^2}\right),$ (ii) $\prod_{p \leq x}(1 - 1/p) < \dfrac{e^{-\gamma}}{\log x}\left(1 + \dfrac{1}{2(\log x)^2}\right),$

(iii) $\prod_{p \leq x} p/(p-1) > e^{\gamma} \log x \left(1 - \dfrac{1}{2(\log x)^2}\right),$ (iv) $\prod_{p \leq x} \dfrac{p}{p-1} < e^{\gamma} \log x \left(1 + \dfrac{1}{2(\log x)^2}\right).$

Proof: These are improved versions of the original works, see [RS] for more details. ∎

**Proposition 36.** Let $x \geq 2$ be a real number, then $\prod_{p \leq x}(1 + 1/p) = \dfrac{6e^{\gamma}}{\pi^2} \log x + O(1/x \log x)$.

Proof: Using standard analysis one has

$$\sum_{p \le x} \log\left(1 + \frac{1}{p}\right) = \sum_{p \le x}\sum_{n=1}^{\infty} \frac{(-1)^{n+1}}{np^n} = \sum_{p \le x}\sum_{n=2}^{\infty} \frac{(-1)^{n+1}}{np^n} + \sum_{p \le x}\frac{1}{p}$$

$$= \sum_{p \le x}\sum_{n=2}^{\infty} \frac{(-1)^{n+1}}{np^n} + \gamma - \sum_{p}\sum_{n=2}^{\infty}\frac{1}{np^n} + \log\log x + R(x)$$

$$= -\sum_{p}\sum_{n=1}^{\infty}\frac{1}{np^{2n}} + \gamma + \log\log x - \sum_{p > x}\sum_{n=2}^{\infty}\frac{(-1)^{n+1}}{np^n} + R(x)$$

$$= \log(6/\pi^2) + \gamma + \log\log x + \sum_{p > x}\sum_{n=2}^{\infty}\frac{(-1)^n}{np^n} + R(x),$$

(7)

where the term $\gamma - \sum_{p}\sum_{n=2}^{\infty}\frac{1}{np^n}$ in the second line is Mertens constant, see Section 15. To complete the proof observe that the error term $E(x) = \sum_{p > x}\sum_{n=2}^{\infty}\frac{(-1)^n}{np^n} + R(x) = O(1/x\log x)$, see Proposition 34.     ∎

The asymptotic version of this result is Theorem 25, see also [CU, p. 110], and [EL, p. 31] for similar details.

## 10 An Estimate of the Totient Function

The key new idea the analysis of the arithmetic functions $N/\varphi(N)$ and $\sigma(N)/N$ is the consideration of the sign and the order of the logarithmic difference as a function of $N$. Previous works in the literature seem to have assumed that the logarithmic difference $\log\log p_k - \log\log\log N = 0$. This assumption appears to places an obstacle in both numerical calculations and in the theory of these functions since any term in this analysis is important. The proof below employs a reductio ad absurdum argument.

**Theorem 7** Let $N_k = 2 \cdot 3 \cdot 5 \cdots p_k$. Then $N_k / \varphi(N_k) > e^\gamma \log\log N_k$ for all sufficiently large integer $N_k$.

Proof: On the contrary suppose that $N_k / \varphi(N_k) \le e^\gamma \log\log N_k$. Then

$$\log \prod_{p \mid N_k}\left(1 - \frac{1}{p^2}\right)^{-1}\left(1 + \frac{1}{p}\right) \le \log\left(e^\gamma \log\log N_k\right),$$

(8)

see Proposition 8-i. Expanding the left side into power series and some analysis yield:

$$\sum_{p \mid N_k}\log\left(1 - \frac{1}{p}\right)^{-1} = \sum_{p \le p_k}\sum_{n=1}^{\infty}\frac{1}{np^n} = \sum_{p \le p_k}\sum_{n=2}^{\infty}\frac{1}{np^n} + \sum_{p \le p_k}\frac{1}{p}$$

$$= \sum_{p \le p_k}\sum_{n=2}^{\infty}\frac{1}{np^n} + \gamma - \sum_{p \ge 2}\sum_{n=2}^{\infty}\frac{1}{np^n} + \log\log p_k + R(N_k)$$

(9)

$$= -\sum_{p > p_k}\sum_{n=2}^{\infty}\frac{1}{np^n} + \gamma + \log\log p_k + R(N_k),$$

where $B = \gamma - \sum_{p \ge 2}\sum_{n=2}^{\infty}\frac{1}{np^n}$ is a constant and $R(N_k)$ is the error term in Mertens formula, see Theorem 12-i. Now inequality (8) becomes

$$- \sum_{p > p_k} \sum_{n=2}^{\infty} \frac{1}{np^n} + \log \log p_k + R(N_k) \leq \log \log \log N_k. \tag{10}$$

Next, since the Chebychev's function satisfies $p_k > \vartheta(p_k) = \log N_k$, and $|p_k - \log N_k| \leq c_5 \log N_k / (\log \log N_k)^A$, for some constants $A > 0$, $c_5 > 0$, one can rewrite the logarithmic difference as

$$\log \log p_k - \log \log \log N_k = O\big(1/(\log \log N_k)^{A+1}\big) > 0, \tag{11}$$

unconditionally, see Theorem 19 and Proposition 21 for more details. Further, using the unconditional form of the error term $R(x) = O(1/\log^B x)$, $B > 0$ constant, see (6), one can rewrite the penultimate line as

$$O(1/(\log \log N_k)^{A+1}) + O(1/(\log \log N_k)^B) \leq \sum_{p > p_k} \sum_{n=2}^{\infty} \frac{1}{np^n}. \tag{12}$$

The left side is of order $O(1/(\log \log N_k)^C) > 0$, where $C = \min \{ A + 1, B \} \geq 1$, and the right side is of order $O(1/\log N_k)$. Ergo for all sufficiently large integers $N_k$, this is a contradiction. ∎

The contradiction in inequality (12) remains in effect for any error term $R(x)$ and any $\log \log p_k - \log \log \log N = O\big(1/f(\log N) \log \log N\big)$ with $\log x \leq f(x) \leq x^{1/2}$, see Proposition 21. Accordingly, the result is consistent with the Riemann hypothesis and numerical data.

## 11 The Divisor Function Inequality

This section proposes a weak form of the divisor function inequality $\sigma(N) < e^\gamma N \log \log N$ for large $N \geq 5041$ due to Ramanujan and Robin. It builds on the analysis of the previous sections and previous papers.

In the previous section the logarithmic difference

$$\log \log p_k - \log \log \log N > 0,$$

is positive for $N = 2^{v_1} \cdot 3^{v_2} \Lambda \ p_k^{v_k} \approx 2 \cdot 3 \Lambda \ p_k$ since $p_k > \log N$. This was used to prove the inequality $N/\varphi(N) > e^\gamma N \log \log N$. On the other hand if $N = 2^{v_1} \cdot 3^{v_2} \Lambda \ p_k^{v_k}$ is sufficiently abundant, then $p_k < \log N$. Accordingly the logarithmic difference

$$\log \log p_k - \log \log \log N < 0$$

is negative. This property will be used in the next result.

The worst case $N = 2^{v_1} \cdot 3^{v_2} \Lambda \ p_k^{v_k}$, where $p_i$ is the $i$th prime, $v_i \geq 1$, will be assumed, the arbitrary case $N = q_1^{v_1} \cdot q_2^{v_2} \Lambda \ q_k^{v_k}$, where $q_i$ is prime and $v_i \geq 0$, is handled in the same way mutatis mutandis.

**Theorem 8.** Let $N \in \mathbb{N}$, then $\sigma(N) < e^\gamma N \log \log N$ for all sufficiently large integers $N \geq 1$ such that $P(N) < \big(1 - 1/9 \log \log N\big) \log N$.

Proof: Let $N = 2^{v_1} \cdot 3^{v_2} \Lambda \; p_k^{v_k} > 5041$ be a sufficiently large abundant integer (possibly colossally abundant) such that $p_k < (1 - 1/9 \log \log N) \log N$ and suppose that $\sigma(N) \geq e^\gamma N \log \log N$. Then

$$\log \prod_{p \,|\, N} \left(1 - \frac{1}{p^2}\right)^{-1} \left(1 + \frac{1}{p}\right) = \log \prod_{p \,|\, N} \left(1 - \frac{1}{p}\right)^{-1} > \log\!\left(e^\gamma \log \log N\right), \tag{13}$$

this follows from Proposition 8-iii. Continuing as in (8) to (10) one has

$$-\sum_{p > p_k} \sum_{n=2}^{\infty} \frac{1}{np^n} + \log \log p_k + R(N) > \log \log \log N \,. \tag{14}$$

Replacing the unconditional estimate of $R(N)$, see (6), and applying Proposition 11 with $f(x) = 9\log x$ return

$$\frac{1}{10(\log \log N)^2} + \frac{4}{15(\log \log N)^3} > \log \log \log N - \log \log p_k + \sum_{p > p_k} \sum_{n=2}^{\infty} \frac{1}{np^n}$$
$$> \frac{1}{9(\log \log N)^2} + \sum_{p > p_k} \sum_{n=2}^{\infty} \frac{1}{np^n} \,. \tag{15}$$

Since the power series on the right side has order $O(1/\log N) > 0$, this is a contradiction for any sufficiently large integer $N > 1$. ∎

## 12 Result for Small Powers of 2

The product $\rho(N) = \prod_{p^\alpha \| N} (1 - 1/p^{\alpha+1})$ which appears in the representation $\sigma(N) = \rho(N) N / \varphi(N)$ is a real-value arithmetic function $\rho : \mathbb{N} \to (6/\pi^2, 1)$. For each finite integer $N$, there is a unique real number in the range $6/\pi^2 < \rho(N) < 1$, (its image $\rho(\mathbb{N})$ is probably dense in the interval $[6/\pi^2, 1]$). The real number $\rho(N)$ is a function of the exponents vector $(v_1, v_2, \ldots, v_k)$ of $N = p_1^{v_1} \cdot p_2^{v_2} \Lambda \; p_k^{v_k}$. Squarefree integers are mapped to the lower end of the interval and numbers divisible by high powers of small primes 2, 3, $\ldots$, are mapped to the upper end of the interval. For example, the value $\rho(N)$ of an extremely abundant number increases toward 1 as $N$ increases but remains bounded $\rho(N) < 1$.

**Theorem 9.** If $N = 2^\alpha M$, where $M$ is odd and $2^\alpha \leq (\log \log N)^2$, then $\sigma(N) < e^\gamma N \log \log N$ for all such $N$ but finitely many exceptions.

Proof: On the contrary suppose that $e^\gamma \log \log N \leq \sigma(N)$. Then

$$e^\gamma \log \log N \leq \frac{\sigma(N)}{N} = \frac{N}{\varphi(N)} \rho(N) < e^\gamma \log \log N \left(1 + \frac{1}{2(\log \log N)^2}\right) \rho(N) \,.$$

The inequality in the right side follows from Theorem 35-iv. Rearranging it returns

$$\rho(N)^{-1} < 1 + \frac{1}{2(\log \log N)^2}$$

Now observe that

$$1 + \frac{1}{2^{\alpha+1} - 1} < \prod_{p^\alpha \| N} (1 - 1/p^{\alpha+1})^{-1} = \rho(N)^{-1} < 1 + \frac{1}{2(\log \log N)^2}.$$

But this contradicts the hypothesis $2^\alpha \leq (\log \log N)^2$. Therefore $\sigma(N) < e^\gamma N \log \log N$ for all $N = 2^\alpha M$ such that $2^\alpha \leq (\log \log N)^2$. ∎

For example, any odd integer $N = M$, or any even integer $N = 2M$, $N = 4M > 10^{24}$, and so on. In fact, almost every integer $N = 2^\alpha M$ satisfies the constraint $2^\alpha \leq (\log \log N)^2$, $M \geq 1$ odd. In contrast, sufficiently abundant integers and sufficiently smooth integers do not satisfy this constraint. For example, those numbers listed in the tables in [AE] and [L], and other integers such as $N = n! \geq 16!$, et cetera, do not satisfy the constraint $2^\alpha \leq (\log \log N)^2$. The integers $N = 2^\alpha M$, $M$ odd, and $2^\alpha > (\log \log N)^2$ can be handled using the inequality

$$\frac{\sigma(N)}{N} < \zeta(3) \prod_{p \mid N} \left(1 + \frac{1}{p} + \frac{1}{p^2}\right)$$

in Proposition 12-iv. The analysis is similar to the proof of Theorem 7.

### 13 s-Free Integers and the Divisor Function Inequality
This section provides a the divisor function inequality for large $s$-free integers $N$. It extends the idea of $s$-free integers as in the paper [S].

**Theorem 37.** Let $N \in \mathbb{N}$ be an $s$-free integer. Then $\sigma(N) < (e^\gamma / \zeta(s)) N \log \log N + O\left(1 / (\log \log N)^2\right)$.

Proof: $N = 2^{v_1} \cdot 3^{v_2} \Lambda \ p_k^{v_k}$, and let $s = \max \{ v_1 + 1, \ldots, v_k + 1 \} \geq 2$. Then

$$\frac{\sigma(N)}{N} = \prod_{p^v \| N} \left(\frac{1 - p^{-v}}{1 - p^{-1}}\right) \leq \prod_{p \mid N} \left(\frac{1 - p^{-s}}{1 - p^{-1}}\right). \tag{16}$$

Rewrite this as

$$\frac{\sigma(N)}{N} \leq \zeta(s) \prod_{p > p_k} \left(1 - p^{-s}\right)^{-1} \prod_{p \leq p_k} \left(1 - p^{-1}\right)^{-1}. \tag{17}$$

Taking logarithms on both sides return

$$\log \frac{\sigma(N)}{N} \leq -\log \zeta(s) + \sum_{p > p_k} \sum_{n=1}^{\infty} \frac{1}{n p^{sn}} + \sum_{p \leq p_k} \sum_{n=1}^{\infty} \frac{1}{n p^n}. \tag{18}$$

Separating the last power series into two sums yields

$$\log \frac{\sigma(N)}{N} \leq -\log \zeta(s) + \sum_{p > p_k} \sum_{n=1}^{\infty} \frac{1}{n p^{sn}} + \sum_{p \leq p_k} \sum_{n=2}^{\infty} \frac{1}{n p^n} + \gamma - \sum_{p \geq 2} \sum_{n=1}^{\infty} \frac{1}{n p^n} + \log \log p_k + R(N). \tag{19}$$

Simplifying this expression yields

$$\log \frac{\sigma(N)}{N} \leq -\log \zeta(s) + \gamma + \log \log p_k + E(N),\qquad(20)$$

where

$$E(N) = R(N) - \sum_{p > p_k} \sum_{n=2}^{\infty} \frac{1}{np^n} + \sum_{p > p_k} \sum_{n=1}^{\infty} \frac{1}{np^{sn}} = O\!\left(1/(\log \log N)^2\right).\qquad(21)$$

Now reverse the logarithm function to complete the claim.  ∎

The special case $s = 2$ in Proposition 36 is well known.

The density of primes in an integer is a deciding factor on the size of the divisor function. The values of the normalized divisor function $\sigma(N)/N$ of integers with low density of prime factors are mapped toward the lower end of the interval $(1, e^{\gamma} N \log \log N)$ and the values $\sigma(N)/N$ of integers with high density of prime factors are mapped toward the upper end of the interval.

**_Theorem_ 38.** (Sierpinski 1956) The arithmetic functions $\sigma(N)/N$ and $\varphi(N)/N$ are dense in $(0, 1)$. In particular, for any $x \in (0, 1)$ and $\varepsilon > 0$ there exist infinitely many sufficiently large $N \in \mathbb{N}$ such that

$$\left| \frac{\sigma(N)}{e^{\gamma} N \log \log N} - x \right| < \varepsilon \quad \text{and} \quad \left| \frac{N}{e^{\gamma} \varphi(N) \log \log N} - x \right| < \varepsilon.$$

## 14 Divisor and Totient Functions on Arithmetic Progressions

It has been known (and proven) for quite sometimes that inequality $\sigma(N) < e^{\gamma} N \log \log N$ holds for odd integers, but remain unproven for even integers, see [SE, p.5].

The natural extension of this observation is to derive results as the following:

$$\frac{\sigma(N)}{N} < \frac{c_0(a,q)}{\varphi(q)} e^{\gamma} \log \log N \quad \text{and} \quad \frac{\varphi(N)}{N} > \frac{c_1(a,q)}{\varphi(q)} \frac{1}{e^{\gamma} \log \log N}$$

for every integer $N \equiv a \bmod q$, where $c_i(a,q) > 0$ is a constant.

For $q = 2$, it appears that $c_0(0,q) = 1$, $c_0(1,q) = 1/2$, and $c_1(0,q) = 1/2$, $c_1(1,q) = 1$.

## 15 An applications to Sums of Four Squares Representations

It is not obvious at all which integers have the maximum number of representations as sums of four squares. For example, does $5040 = x^2 + y^2 + z^2 + w^2$ have more solutions $(x, y, z, w)$ than $5041 = r^2 + s^2 + t^2 + u^2$? Furthermore, it is not obvious why an additive problem about an integer $N$ is deeply influenced by the multiplicative structure of $N$.

**Theorem 39.** (Lagrange 1770) Every integer has a representation as a sum of four squares.

**Theorem 40.** (Jacobi 1834) The number of representations of an integer as a sum of four squares is given by

$$r_4(N) = \begin{cases} 8\sigma(N) & \text{if } N \text{ is odd,} \\ 24\sigma(M) & \text{if } N = 2^\alpha M, \ M \text{ is odd and } \alpha > 0. \end{cases}$$

**Theorem 41.** Let $N = 2^\alpha M > 15$, $M$ odd, be an integer. Then the followings hold.
(i) The number of solutions in the Diophantine equation $N = x^2 + y^2 + z^2 + w^2$ is given by

$$r_4(N) < \begin{cases} 4e^\gamma N \log\log N & \alpha = 0, \\ 24e^\gamma M \log\log M & \alpha \neq 0. \end{cases}$$

(ii) The number of solutions of the Diophantine equation is normally distributed with mean $\mu = \omega(N)$ and standard deviation $\sigma = \sqrt{\log\log N}$.

(iii) Almost every integer has $r_4(N) = \begin{cases} O(N \log\log\log N) & \alpha = 0, \\ O(M \log\log\log M) & \alpha \neq 0, \end{cases}$ solutions.

Proof: Statement (i) follows from Theorem 9 and 40. Statement (ii) follows from Theorem 18 and 40. Lastly, Statement (iii) follows from Theorem 26 and 40.                         Quod erat demonstrandum.  ∎

## 16 Useful Formulae
*Power Series*

$$\log(1+x) = \sum_{n=1}^{\alpha} \frac{(-1)^{n+1}}{n} x^n = x - \frac{x^2}{2} + \frac{x^3}{3} - \Lambda \ , \qquad e^x = \sum_{n=1}^{\alpha} \frac{x^n}{n!} = 1 + x + \frac{x^2}{2!} + \frac{x^3}{3!} + \Lambda$$

For every inter $m \geq 0$, the inequalities hold:

$$\log(1+x) \leq \sum_{n=1}^{2m} \frac{(-1)^{n+1} x^n}{n} \quad \text{and} \quad \log(1+x) \geq \sum_{n=1}^{2m+1} \frac{(-1)^{n+1} x^n}{n}$$

*Constants*
The Euler's constant $\gamma = 0.57721\ 56649\ 01532\ 86060\ 65120\ 90082\ 40243\ 10421\ 59335\ 93992\ldots$ is one of the most important constant in the mathematical sciences. This constant is defined by the limit formula $\gamma = \lim_{x\to\infty} \sum_{n\leq x} 1/n - \log x$. There are dozens of integrals, power series, limits, rapidly convergent series, and other forms of representations of this constant. For example, $\gamma = \lim_{x\to\infty} \sum_{p\leq x} \log\left(1 + \frac{1}{p-1}\right) - \log\log x$ is another representation. The next obvious generalization of this constant is to arithmetic progressions. The definition takes the form $\gamma_{a,q} = \lim_{x\to\infty} \sum_{n\leq x, \ n\equiv a \bmod q} 1/n - \frac{1}{\varphi(q)} \log x$, where $a, q$ are constants. There are other generalizations of the Euler's constant, see the literature.

Meissel-Mertens constant is $B = 0.2614972128476427837554268386086958590515 6664...$. The definition of the constant is by the limit $B = \lim_{x \to \infty} \sum_{p \leq x} 1/p - \log \log x$. The two constants are linked in the formula

$B = \gamma + \sum_{p \geq 2} \log(1 - 1/p) + 1/p$, a proof appears in [HW, p. 351]. The original constant was stated using the

rapidly convergent series $B = \gamma + \sum_{n \geq 2} \frac{\mu(n)}{n} \zeta(n)$, see [VL]. The Mertens constant on an arithmetic progression is defined by $B_{a,q} = \lim_{x \to \infty} \sum_{p \leq x, \, p \equiv a \bmod q} 1/p - \frac{1}{\varphi(q)} \log \log x$, where $a$, $q$ are constants, see [FN] and [LZ] for more details.

*Special Values of the Zeta Function*

Some special values of the zeta function $\zeta(s) = \sum_{n \geq 1} 1/n^{-s} = \prod_p (1 - p^{-s})^{-1}$ are recorded here.

$\zeta(2) = \pi^2/6$, $\zeta(4) = \pi^4/90$, $\zeta(6) = \pi^6/945, ...,$

$\zeta(3) = 1.2020569031595942853997381615114499907649862923 4049...,$

$\zeta(5) = 1.0369277551...$, $\zeta(7) = 1.0083492774...$.